\documentclass{amsart}
\usepackage{amsmath,amsfonts}
\usepackage{latexsym}
\usepackage[latin1]{inputenc}
\title[Analytic Operators]{Operators with analytic orbit for the torus action}
\author[R. Cabral, S. Melo]{Rodrigo A. H. M. Cabral \and Severino T. Melo}
\newtheorem{thm}{Theorem}

\newtheorem{cor}{Corollary}

\newcommand{\R}{{\mathbb R}}

\newcommand{\N}{{\mathbb N}}

\newcommand{\Z}{{\mathbb Z}}
\newcommand{\rn}{{\mathbb R}^{n}}
\newcommand{\tn}{{\mathbb T}^{n}}

\newcommand{\Op}{\mbox{Op}}
\newcommand{\lltt}{\mathcal{L}(L^2({\mathbb T}^n))}

\newcommand{\ops}{operators}
\newcommand{\psd}{pseudo\-dif\-fer\-en\-tial}

\newcommand{\cqd}{\hfill$\Box$}

\newcommand{\pf}{{\tt Proof}: }

\begin{document}

\begin{abstract}
  The class of the bounded operators
  on $L^2(\tn)$ with analytic orbit under the action of $\tn$ by conjugation with the translation operators is shown to coincide 
with the class of those pseudodifferential operators on $\tn$ whose discrete symbol $(a_j)_{j\in\Z^n}$ is {\em uniformly analytic}, in the
sense that there exists $C>1$ such that the derivatives of $a_j$ satisfy $|\partial^\alpha a_j(x)|\leq C^{1+|\alpha|}\alpha!$ for all $x\in\tn$,
all $j\in\Z^n$ and all
$\alpha\in\N^n$.
\end{abstract}
\maketitle
\begin{center}{\footnotesize 2010 Mathematics Subject Classification. Primary 47G30; Secondary 35S05, 58G15.}
\end{center}

\section*{Introduction}

Let us consider the unitary representations $y\mapsto T_y$ and $\eta\mapsto M_\eta$ of $\rn$ on $L^2(\rn)$ defined by $T_yu(x)=u(x-y)$ and 
$M_\eta u(x)=e^{ix\cdot\eta}u(x)$, $x\in\R$. Cordes \cite{C} proved that a bounded linear operator $A$ on $L^2(\rn)$ is such that
\begin{equation}
  \label{hsmooth}
  \R^{2n}\ni (y,\eta)\mapsto M_\eta T_yAT_{-y}M_{-\eta}\in\mathcal{L}(L^2(\rn))
\end{equation}
is a smooth function with values in the Banach space  of all 
bounded operators on $L^2(\rn)$ if and only if there is a $p\in C^{\infty}(\R^{2n})$, bounded and with all its partial derivatives also bounded,
such that, for all smooth and rapidly decreasing $u$ and all $x\in\rn$, one has
\begin{equation}\label{000}
Au(x)=\frac{1}{(2\pi)^n}\int_{\rn}e^{ix\cdot\xi}p(x,\xi)\widehat{u}(\xi)\,d\xi,\ \ \ \widehat{u}(\xi)=\int_{\rn} e^{-is\cdot\xi}u(s)\,ds.
\end{equation}
This class of operators is often denoted by $OPS^0_{0,0}(\rn)$, they are the pseudodifferential operators of order zero with 
symbols satisfying $\rho=\delta=0$ H\"ormander-type global estimates.
Cordes' result can thus be regarded as the characterization of $OPS^0_{0,0}(\rn)$ as the operators with smooth orbit under 
the canonical (not everywhere strongly continuous) action of the Heisenberg group on the C$^*$-algebra $\mathcal{L}(L^2(\rn))$.
It implies, in particular, that this class is a spectrally invariant *-subalgebra of $\mathcal{L}(L^2(\rn))$.

The class $OPS^0_{0,0}(\rn)$ is not invariant under diffeomorphisms. Then it does not make sense, in general,
to define it on manifolds. But since $OPS^0_{0,0}(\rn)$ is invariant under translations, it does make sense to define $OPS^0_{0,0}(\tn)$
as the class of operators acting on functions defined on the torus $\tn$ which are ``locally'' 
(the quotation marks indicate that only the canonical charts of the torus are considered)
given by operators in $OPS^0_{0,0}(\rn)$. It follows from  Cordes' result on $\rn$ that $OPS^0_{0,0}(\tn)$ can be characterized 
as those bounded operators on $L^2(\tn)$ which have smooth orbits when acted by the group $\tn$ via conjugation with the translations. 
This is stated, using the discrete symbol representation of pseudodifferential operators on $\tn$, in \cite[Theorem 2]{M} for the case $n=1$ 
and here in our Theorem~\ref{smooth}.

The main purpose of this paper is to address the question of which bounded operators on $L^2(\tn)$ have real analytic orbits. We prove in 
Theorem~\ref{analytic} that they consist precisely of the operators in $OPS^0_{0,0}(\tn)$ whose discrete symbols, which are sequences $(a_j)$, 
are ``uniformly analytic'', meaning that each $a_j$ is analytic and that the coefficients of their Taylor series 
satisfy estimates uniform in $j$. Following \cite{CC}, we denote by $\mathfrak{Op}_{0}(\tn)$ this class of pseudodifferential operators.

The global representation of pseudodifferential operators on the torus, in which the discrete Fourier transform and discrete symbols
replace, respectively, the Fourier transform on $\rn$ and localization, goes back to Volevich in the 1970's \cite{A}. 
A complete treatment of this subject, including analogues of the Fourier transform for noncommutative Lie groups, is given in \cite{RT}. 

It is an easy  corollary of our Theorem~\ref{analytic} that $\mathfrak{Op}_{0}(\tn)$ is a spectrally invariant *-subalgebra of $\lltt$. 
Chinni and Cordaro \cite{CC}\ give a more concrete proof that it is a *-algebra, but we don't know if it there is a proof of the spectral 
invariance directly using local or discrete symbols. It also follows that $\mathfrak{Op}_{0}(\tn)$ is dense in $OPS^0_{0,0}(\tn)$.

Characterizations of classes of pseudodifferential operators as the smooth operators for actions of Lie groups on C$^*$-algebras have
also been considered in \cite{MM,MM2,P,Ta}, in connection with questions arising from Deformation Quantization and Nonlinear PDEs. 

\section{Smoothness result}\label{s1}

Let $\tn$ denote the torus $\rn/(2\pi\Z)^n$. For each $j\in\Z^n$, let $e_j\in C^\infty(\tn)$ be defined by $e_j(x)=e^{ij\cdot x}$
We have just denoted, as we often will, by the same letter $x$ both an element $x$ of $\rn$ and its class $[x]\in\tn$. We equip
$\tn$ with the measure induced by the Lebesgue measure of $\rn$. If $g\in L^1 (\,]\!-\pi,\pi]^n)$ and $f([x])=g(x)$, we then have
    \[
    \int_{\tn}f=\int_{ ]-\pi,\pi]^{n}}g(x)\,dx.
    \]

Given a linear map
$A:C^\infty(\tn)\to C^\infty(\tn)$, the sequence $(a_j)_{j\in\Z^n}$, $a_j=e_{-j}(Ae_j)\in C^\infty(\tn)$, is called the
{\em discrete symbol} of $A$. Using Fourier series, one may then write
\begin{equation}\label{formula}
  Au(x)=\frac{1}{(2\pi)^n}\sum_{j\in\Z^n}a_j(x)e_j(x)\widehat{u}_j,\ \ \ \widehat{u}_j=\int_{\tn}e_{-j}u, 
  \end{equation}
for all $u\in C^\infty(\tn)$ and $x\in\tn$.

We say that $(a_j)_{j\in\Z^n}$, $a_j\in C^\infty(\tn)$, is a {\em symbol of order} $m\in\R$ if, for every multiindice
  $\alpha\in\N^n$, 
  \[
 \sup\{(1+|j|)^{-m}|\partial^\alpha a_j(x)|;\ j\in\Z^n,\ x\in\tn\}<\infty, 
 \]
where $|j|=(j_1^2+\cdots+j_n^2)^{1/2}$, $j=(j_1,\dots,j_n)$. 
For $n=1$, the following theorem is Proposition 1 of \cite{M}. To adapt that proof to the case $n>1$, one uses  
$e^{ij\cdot x}=(1+|j|^2)^{-p}|(1-\Delta)^pe^{ij\cdot x}|,\ j\in\Z^n,\ p\in\N,$
where $\Delta$ denotes the Laplace operator. See also \cite[Theorem 3.4]{CC}. 

 \begin{thm}\label{t1} If $(a_j)_{j\in\Z^n}$ is a symbol of order $m$, then the formula $(\ref{formula})$ defines a
   linear operator $A:C^\infty(\tn)\to C^\infty(\tn)$. If $m\leq 0$,
 then $A$ extends to a bounded linear operator $A:L^2(\tn)\to L^2(\tn)$, whose norm satisfies, for any integer $p>\frac{n}{2}$,
 \begin{equation}
\label{norm}
||A||\leq C_p\sup\{|(1-\Delta)^pa_j(x)|;\,j\in\Z^n,\,x\in\tn\},\ C_p=\sum_{l\in\Z^n}(1+|l|^2)^{-p}.
\end{equation}
We then write $A=\Op((a_j)_{j\in\Z^n})$, or simply $A=\Op(a_j)$.
\end{thm}

For each $y\in\tn$, let $T_y$ denote the operator, unitary on $L^2(\tn)$, defined by $(T_yu)(x)=u(x-y)$. For $A=\Op(a_j)$, one has
$T_yAT_{-y}=\Op((T_ya_j)_{j\in\Z^n})$. For each $\alpha\in\N^n$, let us denote $A^\alpha=\Op((\partial^\alpha a_j)_{j\in\Z^n})$. 
It follows from the group property $T_yT_z=T_{y+z}$ and from the estimate (\ref{norm}) that $y\mapsto T_yAT_{-y}$ is 
a smooth function on $\tn$ with values in the Banach space $\lltt$ and, moreover, 
$\partial_y^\alpha (T_yAT_{-y})=T_yA^\alpha T_{-y}$. This proves the {\em if} statement in the following theorem.

\begin{thm}\label{smooth}
  A bounded operator $A\in\lltt$ is such that the map
  \begin{equation}
    \label{tsmooth}
    \tn\ni y\mapsto T_yAT_{-y}
  \end{equation}
  is smooth with respect to the norm topology of 
$\lltt$ if and only if $A=\Op(a_j)$ for some symbol $(a_j)_{j\in\Z^n}$ of order zero.
\end{thm}

\pf If a bounded operator $A$ on $L^2(\rn)$ is such that the map (\ref{hsmooth}) is smooth, we call it a {\em Heisenberg smooth} operator.
Analogously, if $A\in\lltt$ is such that (\ref{tsmooth}) is smooth, we say that it is {\em translation smooth}.

For each $\theta=(\theta_1,\cdots,\theta_n)\in\rn$, let $\tilde{U}_\theta\subset\rn$ denote the open rectangle
\[
\tilde{U}_\theta=]\theta_1-\pi,\theta_1+\pi[\times\cdots\times]\theta_n-\pi,\theta_n+\pi[,
\]
let $U_\theta\subset\rn/(2\pi\Z)^n=\tn$ denote the set of the equivalence classes of all $x\in\tilde{U}_\theta$ and let
$\chi_\theta:U_\theta\to\tilde{U}_\theta$ be the chart that
sends a class in $U_\theta$ to its unique representative in $\tilde{U}_\theta$. Let then
$I_\theta:L^2(\tn)\to L^2(\tilde{U}_\theta)$ denote the unitary map $I_\theta u=u\circ\chi_\theta^{-1}$,
$E_\theta:L^2(\tilde{U}_\theta)\to L^2(\rn)$ denote the extension-by-zero embedding and $R_\theta:L^2(\rn)\to L^2(\tilde{U}_\theta)$
denote the restriction map. Finally, let $\Xi_\theta:\mathcal{L}(L^2(\tn))\to\mathcal{L}(L^2(\rn))$  denote the isometric *-homomorphism
$
\Xi_\theta(A)= E_\theta I_\theta AI_\theta^{-1}R_\theta.
$

There exists a finite collection of smooth functions $\{\phi_1,\cdots,\phi_q\}$ (we may choose $q=3^n$) such that $\sum_{r=1}^q\phi_r\equiv 1$ and,
for each pair of integers $1\leq m,p\leq q$, there is a $\theta_{mp}\in\rn$ such that the union of the supports of $\phi_m$ and $\phi_p$ is
contained in $U_{\theta_{mp}}$. Given a translation smooth operator $A\in\lltt$, it can thus be written as a finite sum of operators of
the form $\phi(M) A \psi(M)$, where $\phi(M)$ and $\psi(M)$ are the operators of multiplication by, respectively, the smooth functions
$\phi$ and $\psi$, whose supports are both contained in $U_\theta$ for some $\theta\in\rn$. It is then enough to show that such a
$\phi(M) A \psi(M)$ equals $\Op(a_j)$ for some symbol $(a_j)_{j\in\Z^n}$ of order zero.

Noting that $T_y\phi(M) T_{-y}=\phi_y(M)$, where $\phi_y=T_y\phi$, we see that the operators $\phi(M)$ and $\psi(M)$ are translation smooth and,
hence, that $\phi(M) A \psi(M)$ is also translation smooth. Next we show that $\tilde{A}=\Xi_\theta[\phi(M) A \psi(M)]$ is Heisenberg smooth.
For that it suffices to show (see \cite[Section 8.1]{Cbook}) that both maps 
\begin{equation}
  \label{tilde}
\rn\ni y\mapsto T_y\tilde{A}T_{-y}\ \ \ \mbox{and}\ \ \ \rn\ni\eta\mapsto M_\eta \tilde{A}M_{-\eta}
\end{equation}
are smooth (abusing notation, we are denoting by $T_y$ translation operators both on $L^2(\tn)$ and on $L^2(\rn)$). It suffices to prove
that the partial derivatives exist at $y=0$ and at $\zeta=0$, respectively, because $y\mapsto T_y$ and $\zeta\mapsto M_\zeta$ are representations. 

For $i=1,\dots,n$, let $A_i$ denote the partial derivative with respect to $y_i$ at $y=0$ of $y\mapsto T_yAT_{-y}$,
\begin{equation}
  \label{Newton}
A_i=\lim_{h\to 0}\frac{T_{(0,\dots,0,h,0,\dots 0)}\phi(M)A\psi(M)T_{(0,\dots,0,-h,0,\dots,0)}-\phi(M)A\psi(M)}{h}.
\end{equation}
Since
\[
\Xi_\theta[T_y\phi(M)A\psi(M)T_{-y}]=T_y\Xi_\theta[\phi(M)A\psi(M)]T_{-y}
\]
for all sufficiently small $y$, we conclude from (\ref{Newton}) that the derivative with respect to $y_i$ at $y=0$ of
$y\mapsto T_y\tilde{A}T_{-y}$ equals $\Xi_\theta(A_i)$. 

Let $\bar{\phi}$ and $\bar{\psi}$ be smooth functions with support contained in $U_\theta$ such that, for all sufficiently small $y$,
$\bar{\phi}_y\phi\equiv\phi$ and $\bar{\psi}_y\psi\equiv\psi$, and hence 
\[
\bar{\phi}(M)[T_y\phi(M)A\psi(M)T_{-y}]\bar{\psi}(M)=T_y\phi(M)A\psi(M)T_{-y}.
\]
Multiplying both sides of (\ref{Newton}) on the left by $\bar{\phi}(M)$ and on the right by $\bar{\psi}(M)$, we then conclude that
$
\bar{\phi}(M)A_i\bar{\psi}(M)=A_i
$.
The argument we used to find the first order derivatives can next be used to find the second order
derivatives and, successively, to prove that $\rn\ni y\mapsto T_y\tilde{A}T_{-y}$ is $C^\infty$.

Take $\rho\in C_c^\infty(\tilde{U}_\theta)$ such that $\rho(M)\tilde{A}\rho(M)=\tilde{A}$. The second map in (\ref{tilde}) is smooth because
\[
M_\eta \tilde{A}M_{-\eta}=[M_\eta \rho(M)]\,\tilde{A}\,[\rho(M)M_{-\eta}]
\]
and the maps $\eta\mapsto M_\eta \rho(M)$ and $\eta\mapsto\rho(M)M_{-\eta}$ are smooth.

By Cordes's result \cite{C} (see also \cite[Theorem VIII.2.1]{Cbook} and \cite[Theorem 1]{MM2}), there exists
$p\in C^{\infty}(\R^{2n})$, bounded and with all its partial derivatives also bounded, such that (\ref{000}) holds. The discrete
symbol $(a_j)$ of $\phi(M)A\psi(M)$ can then be expressed in terms of $p$ by the iterated integral
\[
a_j(x)=e_{-j}(x)\,\phi(x)\, A(\psi e_j)(x)=\frac{1}{(2\pi)^n}\int_{\rn}p(x,\xi)\int_{\tilde{U}_\theta}e^{i(x-s)\cdot(\xi-j)}\rho(s)\,ds\,d\xi.
\]
Using a standard technique of pseudodifferential calculus (for $n=1$, what we need is the estimate (10) for $l=0$ in \cite{M}), it follows
that $(a_j)$ is a zero-order discrete symbol. \cqd

\section{Analytic operators}\label{s2}

As in the case of complex valued functions (see, for example, \cite[Section II.3]{J} or \cite[Section V.1]{T}), a smooth function
$f:\tn\to\lltt$ is analytic, in the sense that at each point its Taylor series converges to $f$, if and only if there is $C>1$
such that, for every multiindice $\alpha=(\alpha_1,\dots,\alpha_n)\in\N^n$, 
\begin{equation}
  \label{analyticity}
\sup\{||\partial^\alpha f(x)||;\,x\in\tn\}\leq C^{1+|\alpha|}\alpha!\, 
\end{equation}
(as usual, we denote $|\alpha|=\alpha_1+\dots+\alpha_n$ and $\alpha!=\alpha_1!\dots\alpha_n!$). 

We say that $A\in\lltt$ is {\em translation analytic} if the function $f(y)=T_yAT_{-y}$, $y\in\tn$, is analytic. Being translation smooth,
a translation analytic operator is, by Theorem~\ref{smooth},
necessarily of the form $A=\Op(a_j)$, for some symbol $(a_j)_{j\in\Z^n}$ of order zero.
The following theorem gives a necessary and sufficient condition on $(a_j)$ for $A$ to be translation analytic. 

For each multiindice $\beta=(\beta_1,\cdots,\beta_n)$, 
$L^\beta=(1+\partial_1)^{\beta_1}\dots(1+\partial_n)^{\beta_n}$ defines an invertible operator $L^\beta:C^\infty(\tn)\to C^\infty(\tn)$. Its
inverse can be expressed using Fourier series:
\begin{equation}
  \label{inverse}
(L^\beta)^{-1}u(x)=\frac{1}{(2\pi)^n}\sum_{l\in\Z^n}\frac{e^{i l\cdot x}\widehat{u}_l}{(1+il_1)^{\beta_1}\dots(1+il_n)^{\beta_n}}.
\end{equation}

For each positive integer $m$, let us define a seminorm $\rho_m$ on $C^\infty(\tn)$ by
\[
\rho_m(a)=\sup\{|\partial^\beta a(x)|;\,|\beta|\leq m,\,x\in\tn\}.
\]

\begin{thm}\label{analytic} An $A\in\lltt$ is a translation analytic operator if and only if $A=\Op(a_j)$ for some symbol of order zero
  $(a_j)_{j\in\Z^n}$ satisfying, for some constant $C>1$ and for every multiindice $\alpha$,
  \begin{equation}
    \label{CCsymbol}
  \sup\{|\partial^\alpha a_j(x)|;\,x\in\tn,\,j\in\Z^n\}\leq C^{1+|\alpha|}\alpha!\,.
\end{equation}
  \end{thm}

\pf 
Suppose first that $A=\Op(a_j)$ for some $(a_j)$ satisfying (\ref{CCsymbol}) and consider $f(y)=T_yAT_{-y}$, $y\in\tn$.
Using the remarks and the notation of the paragraph right after Theorem~\ref{t1}, and by the norm estimate (\ref{norm}),
we have, for $p>\frac{n}{2}$,
\[
\sup\{||\partial^\alpha f(x)||;\,x\in\tn\}\leq ||A^\alpha||\leq 2^pC_p\sup_{j\in\Z^n}\rho_{2p}(\partial^\alpha a_j).
\]
It then follows from (\ref{CCsymbol}) that
\[
\sup\{||\partial^\alpha f(x)||;\,x\in\tn\}\leq 2^pC_p C^{1+2p+|\alpha|}[(\alpha_1+2p)!\cdots(\alpha_n+2p)!]
\]\[
\leq\mu^n 2^pC_pC^{2p+1} (2C)^{|\alpha|}\alpha!,
\]
where $\mu=\sup_{t>0}2^{-t}(t+2p)(t+2p-1)\dots(t+1)$. This implies $A$ is translation analytic.

Conversely, let $A$ be translation analytic. By Theorem~\ref{smooth}, $A=\Op(a_j)$ for some zero order discrete symbol $(a_j)$.
For each nonzero multiindice $\beta$, define:
\[
B^\beta=\Op[(L^\beta a_j)_{j\in\Z^n}].
\]
Then $T_yB^\beta T_{-y}=L^\beta f(y)$, where $f(y)=T_yAT_{-y}$, $y\in\tn$, is analytic, by hypothesis.
Since $C>1$, we get from (\ref{analyticity}):
\[
||B^\beta||\leq 2^{|\beta|}C^{1+|\beta|}\beta!.
\]

On the other hand, the discrete symbol of $A$ can be recovered from the discrete symbol $(b^\beta_j)_{j\in\Z^n}$ of $B^\beta$,
$b^\beta_j=e_j(B^\beta e_{-j})=L^\beta a_j$, using (\ref{inverse}):
\[
a_j(x)=\frac{1}{(2\pi)^n}  \sum_{l\in\Z^n}\frac{e^{ilx}[e_j(B^\beta e_{-j})]^{^\wedge}_{_l}}{(1+il_1)^{\beta_1}\dots (1+il_n)^{\beta_n}}.
\]

Given now a multiindice $\alpha$, let $\beta=\alpha+(2,\dots,2)$. Using that $|[e_j(B^\beta e_{-j})]^{^\wedge}_{_l}|\leq (2\pi)^n||B^\beta||$ for all
$j$ and $l$, we get
\[
|\partial^\alpha a_j(x)|=\frac{1}{(2\pi)^n}
\left|\sum_{l\in\Z^n}\frac{(il)^\alpha e^{il\cdot x} [e_j(B^\beta e_{-j})]^{^\wedge}_{_l}}{(1+il_1)^{\beta_1}\dots (1+il_n)^{\beta_n}}\right|
\]
\[
\leq ||B^\beta||\sum_{l\in\Z^n}\frac{|l_1|^{\alpha_1}\dots |l_n|^{\alpha_n}}{(1+l_1^2)^{\beta_1/2}\dots(1+l_n^2)^{\beta_n/2}}
\]
\[
\leq \left(\sum_{p=-\infty}^\infty \frac{1}{1+p^2}\right)^n2^{|\alpha|+2n}C^{1+2n+|\alpha|}\beta!.
\]
We may now bound $\beta!$ by a constant times $\alpha!$, similarly as we did in the first part of this proof, which will then prove (\ref{CCsymbol}).
\cqd

As in \cite{CC}, we denote by $\mathfrak{Op}_{0}(\tn)$ the class of all $A=\Op(a_j)$ with $(a_j)_{j\in\Z^n}$ satisfying (\ref{CCsymbol}).

\begin{cor}\label{corolario}
  The class $\mathfrak{Op}_{0}(\tn)$ is *-subalgebra of $\lltt$ which is spectrally invariant, in the sense  that, if an
  $A\in\mathfrak{Op}_{0}(\tn)$ is invertible in $\lltt$, then its inverse $A^{-1}$ belongs to $\mathfrak{Op}_{0}(\tn)$.
\end{cor}

\pf If $f:\tn\to\lltt$ and $g:\tn\to\lltt$ are real analytic, their pointwise product is analytic. It follows that the product of two
translation analytic operators is translation analytic. The map $y\mapsto f(y)^*$ is also analytic. Hence, the adjoint of a
translation analytic operator is translation analytic. Moreover if $f(y)$ is invertible for every $y\in\tn$,
then $y\mapsto f(y)^{-1}$ is also analytic (because $T_{-y}=T_y^{-1}$ and a pointwise invertible analytic function has 
an analytic inverse). Then, if $A\in\lltt$ is invertible and translation analytic, so is $A^{-1}$.
\cqd

\begin{cor} $\mathfrak{Op}_{0}(\tn)$ is $L^2$-operator-norm dense in  $OPS^0_{0,0}(\tn)$.
\end{cor}
\pf It is the content of \cite[Theorem 4]{N} that the set of analytic vectors of a strongly continuous representation of a Lie group on a
Banach space $\mathfrak{X}$ is dense in $\mathfrak{X}$. We apply that theorem to the representation $y\mapsto T_yAT_{-y}$ of $\tn$  on
the Banach space $\mathfrak{X}=\{A\in\lltt;\,y\mapsto T_yAT_{-y}$ is continuous$\}$. By our Theorem~\ref{analytic}, the set of analytic vectors for
this representation is equal to $\mathfrak{Op}_{0}(\tn)$, which is dense in $\mathfrak{X}$ by Nelson's result. On the other hand, $\mathfrak{X}$
contains the set of smooth vectors, which is equal to $OPS^0_{0,0}(\tn)$ by our Theorem~\ref{smooth}. \cqd

\vskip0.5cm

\tiny{
\noindent
Instituto de  Matem\'atica e Estat\'{\i}stica, Universidade de S\~ao Paulo, Rua do Mat\~ao 1010, 05508-090 S\~ao Paulo, Brazil.
E-mails: rahmc@ime.usp.br, toscano@ime.usp.br.}

\end{document}